\documentclass[12pt]{amsart}

\input xy
\xyoption{all}

\usepackage{geometry}
\usepackage{amsmath}
\usepackage{amssymb}
\usepackage{hyperref}
\usepackage{cite}
\usepackage[final]{showkeys}
\usepackage{hyperref}
\usepackage{setspace}
\usepackage{amsthm}  
\newtheorem{same}{This should never appear}[section]
\newtheorem{defin}[same]{Definition}
\newtheorem{claim}[same]{Claim}

\newtheorem{theorem}[same]{Theorem}
\newtheorem{example}[same]{Example}
\newtheorem{lemma}[same]{Lemma}

\newtheorem{cor}[same]{Corollary}
\newtheorem{prop}[same]{Proposition}

\newbox\noforkbox \newdimen\forklinewidth
\setbox0\hbox{$\textstyle\smile$}
\setbox1\hbox to \wd0{\hfil\vrule width \forklinewidth depth-2pt
 height 10pt \hfil}
\setbox\noforkbox\hbox{\lower 2pt\box1\lower 2pt\box0\relax}
\def\unionstick{\mathop{\copy\noforkbox}\limits}

\def\nonfork_#1{\unionstick_{\textstyle #1}}

\setbox0\hbox{$\textstyle\smile$}
\setbox1\hbox to \wd0{\hfil{\sl /\/}\hfil}
\setbox2\hbox to \wd0{\hfil\vrule height 10pt depth -2pt width
              \forklinewidth\hfil}
\newbox\doesforkbox
\setbox\doesforkbox\hbox{\lower 2pt\box1 \lower 2pt\box2\lower2pt\box0\relax}
\def\nunionstick{\mathop{\copy\doesforkbox}\limits}

\def\fork_#1{\nunionstick_{\textstyle #1}}

\newcommand{\cFml}{\textrm{cFml }}
\newcommand{\sea}{\mathfrak{C}}

\newcommand{\ba}{\bold{a}}
\newcommand{\bb}{\bold{b}}

\newcommand{\bd}{\bold{d}}

\newcommand{\bx}{\bold{x}}
\newcommand{\by}{\bold{y}}
\newcommand{\bz}{\bold{z}}

\newcommand{\rest}{\upharpoonright}

\newcommand{\qp}{\mathbb{Q}'}
\newcommand{\monus}{\dot{-}}

\newcommand{\formula}{\chi}

\newcommand{\term}{\mathfrak{t}}
\newcommand{\Ll}{\mathbb{L}}

\newcommand{\Mod}{\te{Mod }}

\newcommand{\dense}{\text{Dense}}

\newcommand{\seq}[1]{\langle #1 \rangle}
\newcommand{\te}[1]{\textrm{#1}}

\title[$\dense_\tau$ and $K_{dense}$]{A presentation theorem for continuous logic and Metric Abstract Elementary Classes}
\author{Will Boney}
\email{wboney@math.harvard.edu}
\address{Department of Mathematics, Harvard University, Cambridge, MA, USA}
\date{\today
This work was begun while the author was working towards his PhD under Rami Grossberg and
he is grateful for his guidance and support. Part of this material is based upon work done while
the author was supported by the National Science Foundation under Grant No. DMS-1402191.} 

\begin{document}

\maketitle

\begin{abstract}
In recent years, model theory has widened its scope to include metric structures by considering real-valued models whose underlying set is a complete metric space.  We show that it is possible to carry out this work by giving presentation theorems that translate the two main frameworks (continuous first order logic and Metric Abstract Elementary Classes) into discrete settings (a nice fragment of $\mathbb{L}_{\omega_1, \omega}$ and Abstract Elementary Classes, respectively).  We also translate various notions of classification theory.
\end{abstract}

\section{Introduction}

In the spirit of Chang and Shelah's presentation results (from \cite{chang68} and \cite{sh88}, respectively), we prove a presentation theorem for classes of continuous structures in terms of a class of discrete structures.  This presentation is motivated by the fact that continuous functions are determined by their values on a dense subset of their domain.  Focusing on dense subsets is key because it allows us to drop the requirement that structures be complete, which is not a property expressible by discrete (classical) logic, even in the broader contexts of $\Ll_{\lambda, \omega}$ or Abstract Elementary Classes.

The specific statements of the presentation theorems appear below (see Theorem \ref{goal} for continuous first-order logic and Theorem \ref{Kdense} for Metric Abstract Elementary Classes), but the general idea is the same in both cases: given a continuous language $\tau$, we define a discrete language $\tau^+$ that allows us to approximate the values of the functions and relations by a countable dense subset of values, namely $\mathbb{Q} \cap [0,1]$.  Note that the specification that this dense set (and its completion) is standard already requires an $\Ll_{\omega_1, \omega}(\tau^+)$ sentence, even if we are working in continuous first-order logic.  Then, given a continuous $\tau$-structure $M$ and a nicely dense (see Definition \ref{nicelydensedef} below) subset of it $A$, we can form a discrete $\tau^+$-structure $M_A$ with universe $A$ that encodes all of $M$.


This defines a functor between the continuous class and their discrete approximations that witnesses that these two classes are equivalent.  We investigate this equivalence further by considering various model-theoretic properties and finding analogues for them in the class of approximations.  We consider types, saturation, (and for Metric Abstract Elementary Classes), amalgamation, joint embedding, and ${\bf d}$-tameness.

We present two applications of this functor.  In the realm of Metric Abstract Elementary Classes, it allows us to reduce many foundational questions to the same questions about Abstract Elementary Classes, where there are known answers.  For instance, Theorem \ref{maechanf-thm} answers an open question about the existence of Hanf numbers for Metric Abstract Elementary Classes, and shows that it is the same as for normal Abstract Elementary Classes.  On the discrete side, the functor allows for the applications of results about countable fragments of $\mathbb{L}_{\omega_1, \omega}$ to continuous first-order logic.

Throughout this paper, we assume that the reader is familiar with the basics of the continuous contexts--either first-order or Metric Abstract Elementary Classes--but try to provide references and reminders when discussing the concepts (especially to \cite{fourguys} for the continuous first-order context).

Many of the arguments in the two cases are similar.  In order to avoid repeating  the same arguments twice in slightly different contexts, we provide the details only once.  We have chosen to provide the details for continuous first-order logic because this context often allows more specific formulations of the correspondence.

Dense sets are not quite the right context because they need not be substructures of the larger structure.  Instead, we introduce nicely dense sets to require them to be closed under functions.

\begin{defin}\label{nicelydensedef}
Given a continuous model $M$ and a set $A \subset |M|$, we say that $A$ is \emph{nicely dense} iff $A$ is dense in the metric structure $(|M|, d^M)$ and $A$ is closed under the functions of $M$.  
\end{defin}

A few notational points:
\begin{itemize}

	\item We will want to prove similar results for both ``greater than'' and ``less than.''  In order to avoid writing everything twice, we often use $\square$ to stand in for both $\geq$ and $\leq$.  Thus, asserting a statement for ``$r \square s$'' means that that statement is true both for ``$r \geq s$'' and for ``$r \leq s$.''

	\item Our goal is to translate the real-valued formulas of $\tau$ into classic, true/false formulas of $\tau^+$.  We do this by encoding relations into $\tau^+$ that are intended to specify the value of $\phi$ by deciding if it is above or below each possible value.  To ensure that the size of the language doesn't grow, we take advantage of the separability or $\mathbb{R}$ and only compare each $\phi$ to the rationals in $[0,1]$.  Set $\qp := [0,1] \cap \mathbb{Q}$.
	
	\item For $\tau$-structure, we will use roman letters $M, N$, etc. $\overline{\mathcal{A}}$ for $\tau$-structures coming from our functor.  Typically, we will not distinguish them from their universe, but will write, e. g., $|\overline{\mathcal{A}}|$ for this if necessary for clarity.  For $\tau^+$-structures, these will be denoted $M_A$ (if they come from the functor applied to $(A, M)$) and using script letters $\mathcal{A}, \mathcal{B}$, etc. along with the roman letter $A, B$ to denote the universe.
	
\end{itemize}

We would like to thank Pedro Zambrano and the referee for helpful comments on this paper.  After posting drafts of this paper, the author discovered a similar project by Ackerman \cite{nate-paper} (done independently).  Ackerman also encodes metric notions in $\mathbb{L}_{\omega_1, \omega}$, but does so by encoding the category of complete metric spaces and uniformly continuous maps and having the encoding of metric structures fall out of this, rather than focusing on continuous logic from the start (as we do here).  Additionally, our focus is primarily on applications dealing with classification theory (saturation, etc.), while Ackerman applies his methods to explore set-theoretic absoluteness.

\section{Models and Theories} \label{main-sect}

The main thesis of the presentation of continuous first-order logic is that model-theoretic properties of continuous first order structures can be translated to model-theoretic (but typically quantifier free) properties of discrete structures  that model a specific theory in an expanded language.  We use $\cFml \tau$ to denote the continuous formulas of the language $\tau$.  The main theorem about this presentation is the following:
\begin{theorem}\label{goal}
Let $\tau$ be a continuous language.  Then there is
\begin{enumerate}
	\item[(a)] a discrete language $\tau^+$;
	\item[(b)] an $\Ll_{\omega_1, \omega}(\tau^+)$ theory $\dense_\tau$;
	\item[(c)] a map that takes continuous $\tau$-structures $M$ and nicely dense subsets $A$ to discrete $\tau^+$-structures $M_A$ that model $\dense_\tau$;
	\item[(d)] a map that takes discrete $\tau^+$ structures $\mathcal{A}$ that model $\dense_\tau$ to continuous $\tau$-structures $\overline{\mathcal{A}}$
\end{enumerate}
with the properties that
\begin{enumerate}
	\item $M_A \vDash \dense_\tau$ has universe $A$ and, for any $\ba \in A$, $\phi(\bx) \in \cFml \tau$, $r \in \qp$, and $\square$ standing for $\geq$ and $\leq$, we have
	$$M_A \vDash R_{\phi \square r}[\ba] \iff \phi^M(\ba) \square r$$
	
	\item $A$ is a dense subset of $\overline{\mathcal{A}}$ and, for any $\ba \in A$, $\phi(\bx) \in \cFml \tau$, $r \in \qp$, and $\square$ standing for $\geq$ and $\leq$, we have
	$$\mathcal{A} \vDash R_{\phi \square r}[\ba] \iff \phi^{\overline{\mathcal{A}}}(\ba) \square r$$

	\item these maps are (essentially) each other's inverse.  That is, given any nicely dense $A \subset M$, we have $M \cong_A \overline{M_A}$ and, given any $\tau^+$-structure $\mathcal{A} \vDash \dense_\tau$, we have $\left(\overline{\mathcal{A}}\right)_A = \mathcal{A}$.	
\end{enumerate}
\end{theorem}

We say ``essentially'' in the last clause because completions are not technically unique as the objects selected as limits can vary, but this fairly pedantic point is the only obstacle.

Restricting to dense subsets and their completions has already been considered in continuous first-order logic, where it goes by the name \emph{prestructure} (see \cite[Section 3]{fourguys}).  The key difference here is that, while prestructures are still continuous objects with uniformly continuous functions and relations, $M_A$ is a discrete object with the relations on it being either true or false. 

{\bf Proof:} Our proof is long but straightforward.  First, we will define $\tau^+$ and $\dense_\tau$.  Then, we will introduce the map $(M, A) \mapsto M_A$ and prove it satisfies (1).  After this, we will introduce the other map $\mathcal{A} \mapsto \overline{\mathcal{A}}$ and prove (2).  Finally, we will prove that they satisfy (3).\\

{\bf Defining the new language and theory}\\

We define the language $\tau^+$ to be
$$\seq{F_i^+, R_{\phi(\bx) \geq r}, R_{\phi(\bx) \leq r}}_{i < n_F, \phi(\bx) \in \cFml \tau, r \in \qp}$$
with the arity of $F_i^+$ matching the arity of $F_i$ and the arity of $R_{\phi(\bx)\square r}$ matching $\ell(\bx)$.  Since we only use a full (that is, dense) set of connectives (see \cite[Definition 6.1]{fourguys}), we have ensured that $|\tau^+| = |\tau| + \aleph_0$.

We define $\dense_\tau \subset \Ll_{\omega_1, \omega}(\tau^+)$ to be the universal closure of all of the following formulas ranging over all continuous formulas $\phi(\bz)$ and $\psi(\bz')$, all terms $\term(z, \bz'')$, and all $r, s \in \qp$ and $t \in \qp - \{0\}$.  We have divided them into headings so that their meaning is (hopefully) more clear.  When we refer to specific sentences of $\dense_\tau$ later, we reference the ordering in this list.  As always, a $\square$ in a formula means that it should be included with both a `$\geq$' and a `$\leq$' replacing the $\square$.
\begin{enumerate}

	\item {\bf The ordered structure of $\mathbb{R}$}
	\begin{enumerate}
	
			\item $\neg R_{\phi(\bz) \geq r}(\bx) \to R_{\phi(\bz) \leq r}(\bx)$
			
			\item $\neg R_{\phi(\bz) \leq r}(\bx) \to R_{\phi(\bz) \geq r}(\bx)$
			
			\item If $r > s$, then include $\neg R_{\phi(\bz) \geq r}(\bx) \vee \neg R_{\phi(\bz) \leq s}(\bx)$ \label{order1}
			
			\item If $r \geq s$, then include \label{order2}
			\begin{itemize}
				\item $R_{\phi(\bz) \leq s}(\bx) \to R_{\phi(\bz) \leq r}(\bx)$; and
				\item $R_{\phi(\bz) \geq r}(\bx) \to R_{\phi(\bz) \geq s}(\bx)$
			\end{itemize}
			
			\item $R_{\phi(\bz) \geq r}(\bx) \vee R_{\phi(\bz) \leq r}(\bx)$
			
			\item $\wedge_{n<\omega} \vee_{r, s \in \qp, |r-s|<\frac{1}{n}} R_{\phi \leq r}(\bx) \wedge R_{\phi \geq s}(\bx)$ \label{supinf}

			\item $(\wedge_{n < \omega} R_{\phi(\bz) \geq r - \frac{1}{n}}(\bx)) \to R_{\phi(\bz) \geq r}(\bx)$ \label{stdleft}
			\item $(\wedge_{n < \omega} R_{\phi(\bz) \leq r + \frac{1}{n}}(\bx)) \to R_{\phi(\bz) \leq r}(\bx)$ \label{stdright}

	\end{enumerate}
		
	\item {\bf Construction of formulas}
	\begin{enumerate}
	
		\item $R_{\phi(\bz) \geq 0}(\bx) \wedge R_{\phi(\bz) \leq 1}(\bx)$
		\item $\neg R_{0 \geq t}(\bx) \wedge \neg R_{1 \leq 1 -t}(\bx)$;
		
		\item $R_{\frac{\phi(\bz)}{2} \geq r}(\bx) \leftrightarrow R_{\phi(\bz) \geq 2 r}(\bx)$
		\item $R_{\frac{\phi(\bz)}{2} \leq r}(\bx) \leftrightarrow R_{\phi(\bz) \leq 2r}(\bx)$;
		
		\item $R_{\phi(\bz) \monus \psi(\bz')\geq r}(\bx, \bx') \leftrightarrow \vee_{s \in \qp} (R_{\psi(\bz') \leq s}(\bx') \wedge R_{\phi(\bz) \geq r + s}(\bx))$ \label{gmonus}
		\item $R_{\phi(\bz) \monus \psi(\bz') \leq r}(\bx, \bx') \leftrightarrow \vee_{s \in \qp} (\neg R_{\psi(\bz') \leq s}(\bx') \wedge R_{\phi(\bz)\leq r + s}(\bx))$;
		
		\item $R_{\sup_y \term (y, \by) \leq r}(\bx) \leftrightarrow \forall x R_{\term(y, \by) \leq r}(x, \bx)$ \label{lsup}
		\item $R_{\sup_y \term (y, \by) \geq r}(\bx) \leftrightarrow \wedge_{n < \omega} \exists x R_{\term(y, \by) \geq r - \frac{1}{n}}(x, \bx)$ \label{gsup};
		
		\item $R_{\inf_y \term(y, \by) \leq r}(\bx) \leftrightarrow \wedge_{n < \omega} \exists x R_{\term(y, \by) \leq r + \frac{1}{n}}(x, \bx)$ \label{linf};
		\item $R_{\inf_y \term(y, \by) \geq r}(\bx) \leftrightarrow \forall x R_{\phi(y, \by) \geq r}(x, \bx)$ \label{ginf};
		
		\item $R_{\phi(y, \by) \square r}(\term(\bx'), \bx) \leftrightarrow R_{\phi(\term(\by'), \bx) \square r}(\bx', \bx)$ \label{termred}

	\end{enumerate}
		\item {\bf Metric structure}
		\begin{enumerate}
		
		\item $R_{d(y, y') \leq 0}(x, x') \leftrightarrow x = x'$; \label{metric1}
		
		\item $R_{d(y, y') \square r}(x, x') \leftrightarrow R_{d(y, y') \square r}(x', x)$; \label{metric2}
		
		\item $\wedge_{r \in \qp}( R_{d(y, y') \geq r}(x, x') \to \forall x'' \vee_{s \in \qp \cap [0, r]} R_{d(y, y') \geq s}(x, x'') \wedge R_{d(y, y') \geq r-s}(x'', x'))$ \label{metric3}

	\end{enumerate}
	\item {\bf Uniform Continuity}
	\begin{enumerate}
	
		\item For each $r, s \in \qp$ and $i < \tau_F$ such that $s < \Delta_{F_i}(r)$, we include the sentence
		$$\wedge_{i < n} R_{d(z, z') \leq s}(x_i, y_i) \to R_{d(z, z') \leq r}(F_i(\bx), F_i(\by))$$
		
		\item For each $r, s \in \qp$ and $j < \tau_R$ such that $s < \Delta_{R_j}(r)$, we include the sentence
		$$\wedge_{i < n} R_{d(z, z') \leq s}(x_i, y_i) \to (R_{R_j(\bz) \monus R_j(\bz') \leq r}(\bx, \by) \wedge R_{R_j(\bz) \monus R_j(\bz') \leq r}(\by, \bx)) $$

	\end{enumerate}

\end{enumerate}

We have been careful about the specific enumeration of these axioms because, if the original continuous language is countable, then $\dense_\tau$ is countable.  In particular, we could take the conjunction of it and make it a single $\Ll_{\omega_1, \omega}(\tau^+)$ sentence.  This means that it is expressible in a countable fragment of $\Ll_{\omega_1, \omega}(\tau^+)$.  In general, $\dense_\tau$ is expressible in a $|\tau| + \aleph_0$ sized fragment of $\Ll_{\omega_1, \omega}(\tau^+)$.

Countable fragments are the most well-studied infinitary languages and many of the results in, say, Keisler \cite{keislerw1} use these fragments.  These results are then applicable to countable theories in continuous first-order logic.  Such an application is an improvement on finding generalized indiscernibles: in \cite[Lemma 3.35]{by-uncountcat}, Ben-Yaacov appeals to (essentially) Morley's omitting types theorem to find appropriate indiscernibles that look like a long enough sequence.  Even for countable theories $T$, ``long enough'' is bounded by $\beth_{(2^{\omega})^+}$.  However, using our correspondence (including the analysis of types from Section \ref{type-sec}) and the Hanf number for $\mathbb{L}_{\omega_1, \omega}$ \cite[Theorem 21]{keislerw1}, this bound can be brought down to $\beth_{\omega_1}$.\\

{\bf From continuous to discrete...}\\

This is the easier directions.  We define the structure $M_A$ so that all of the ``intended'' correspondences hold and everything works out well.\\

Suppose we have a continuous $\tau$-structure $M$ and a nicely dense subset $A$.  Now we define an $\tau^+$-structure $M_A$ by
\begin{enumerate}
	\item the universe of $M_A$ is $A$;
	\item $(F^+_i)^{M_A} = F_i^M \rest A$ for $i < n_F$; and
	\item for $r \in \qp$ and $\phi(\bx) \in \cFml \tau$, set
	$$R_{\phi \square r}^{M_A} = \{ \ba \in A : \phi^M(\ba) \square r \}$$
\end{enumerate}

This is an $\tau^+$-structure since $A$ is closed under functions.  The real meat of this part is the following claim, which is (1) from the theorem.

\begin{claim}
 $M_A \vDash \dense_\tau$ and, for any $\ba \in A$, $\phi(\bx) \in \cFml \tau$, $r \in \qp$, and $\square \in \{\geq, \leq\}$, we have
	$$M_A \vDash R_{\phi \square r}[\ba] \iff \phi^M(\ba) \square r$$
\end{claim}
{\bf Proof:}  This is all straightforward.  From the definition, we know that, for any $\ba \in A$ and formula $\phi(\bx) \in \cFml \tau$ and $\square \in \{\geq, \leq\}$, we have
$$M_A \vDash R_{\phi \square r}[\ba] \iff \phi^M(\ba) \square r$$
This gives an easy proof of the fact that $M_A \vDash \dense_\tau$ because they are all just true facts if `$R_{\phi \square r}(\ba)$' is replaced by `$\phi(\ba) \square r$.' \hfill $\dag_{Claim}$\\

{\bf ...and back again}\\

This is the harder direction.  We want to `read out' the $\tau$-structure that $A$ is a dense subset of from the $\tau^+$-structure $\mathcal{A}$.  First, we use the axioms of $\dense_\tau$ to show that we can read out the metric and relations of $\tau$ from the relations of $\tau^+$ and that these are well-defined.  Then we complete $A$ and use the uniform continuity of the derived relations to expand them to the whole structure.  In the first direction, $\dense_\tau$ could have been any collection of true sentences about continuous structures and the real line, but this direction makes it clear that the axioms chosen are sufficient.\\

Suppose that we have an $\tau^+$-structure $\mathcal{A}$ that models $\dense_\tau$.  The following claim is an important step in reading out the relations of the completion of $A$ from $\mathcal{A}$.

\begin{claim}\label{sup=inf}
For any $\phi(\bx) \in \cFml \tau$ and $\ba \in A$, we have
	$$\sup \{ t \in \qp : \mathcal{A} \models R_{\phi(\bx) \leq t}(\ba) \} = \inf \{t \in \qp: \mathcal{A} \models R_{\phi(\bx) \geq t}(\ba) \}$$
\end{claim}

{\bf Proof:} We show this equality by showing two inequalities.

\begin{itemize}
		\item Let $r \in \{ t \in \qp : \mathcal{A} \models R_{\phi(\bx) \leq t}(\ba) \}$ and $s \in \{t \in \qp: \mathcal{A} \models R_{\phi(\bx) \geq t}(\ba) \}$.  Then
		$$\mathcal{A} \models R_{\phi(\bx) \geq r}(\ba) \wedge R_{\phi(\bx) \leq s}(\ba)$$
		Then, since $M^+$ satisfies (\ref{order1}), we must have $r \leq s$.  Thus 
		$$\sup \{ t \in \qp : \mathcal{A} \models R_{\phi(\bx) \leq t}(\ba) \} \leq \inf \{t \in \qp: \mathcal{A} \models R_{\phi(\bx) \geq t}(\ba) \}$$
		\item By (\ref{supinf}), we have
		$$\mathcal{A} \models \wedge_{n < \omega} \vee_{r, s \in \qp; |r-s| < \frac{1}{n}} R_{\phi(\bx) \leq r}(\ba) \wedge R_{\phi(\bx) \geq s}(\ba)$$
		Let $\epsilon > 0$.  Then there is $n_0 < \omega$ such that $\epsilon > \frac{1}{n_0}$.  By the above, there are $r, s \in \qp$ such that $|r-s| < \frac{1}{n_0}$ and
		$$\mathcal{A} \models R_{\phi(\bx) \leq r}(\ba) \wedge R_{\phi(\bx) \geq s}(\ba)$$
		As above, (\ref{order1}) implies $r \geq s$, so we have $r - s < \frac{1}{n_0} < \epsilon$.  Thus $r < s + \epsilon$ and $s \in \{ t \in \qp : \mathcal{A} \models R_{\phi(\bx) \leq t}(\ba) \}$ and $r \in \{t \in \qp: \mathcal{A} \models R_{\phi(\bx) \geq t}(\ba) \}$.  Then, $\inf \{t \in \qp: \mathcal{A} \models R_{\phi(\bx) \geq t}(\ba) \} \leq \sup  \{ t \in \qp : \mathcal{A} \models R_{\phi(\bx) \leq t}(\ba) \} $. \hfill $\dag_{Claim}$\\
	\end{itemize}

The first relation that we need is the metric.  Given $a, b \in A$, we set
\begin{eqnarray*}
D(a, b) &: =& \sup \{ r \in \qp : \mathcal{A} \vDash R_{d(x,y) \geq r}[a, b] \} \\
&=& \inf \{ r \in \qp : \mathcal{A} \vDash R_{d(x, y) \leq r}[a, b] \}
\end{eqnarray*}

These definitions are equivalent by Claim \ref{sup=inf}.  We show that this is indeed a metric on $A$.

\begin{claim}\label{Dmetric}
$(A, D)$ is a metric space.
\end{claim}
{\bf Proof:} We go through the metric space axioms.  Let $a, b \in A$.
	\begin{enumerate}
	
		\item \begin{eqnarray*}
		D(a, b) = 0 &\implies& \inf \{ r \in \qp : \mathcal{A} \vDash R_{d(x, y) \leq r} (a, b) \} = 0 \\
		&\implies& \forall n < \omega, \exists r_n \in \qp \te{ so }\mathcal{A} \models R_{d(x, y) \leq r_n}(a, b) \te{ and }r_n \leq \frac{1}{n}\\
		&\implies_{(\ref{order2})}& \forall n < \omega, \mathcal{A} \models R_{d(x, y) \leq \frac{1}{n}}(a, b)\\
		&\implies_{(\ref{stdright})}& \mathcal{A} \models R_{d(x, y) \leq 0}(a, b)\\
		&\implies& a = b
		\end{eqnarray*}
		\begin{eqnarray*}
		a=b &\implies& \mathcal{A} \models R_{d(x, y) \leq 0}(a, b)\\
		&\implies& \inf \{r \in \qp : \mathcal{A} \models R_{d(x, y) \leq r}(a, b) \} = 0\\
		&\implies& D(a, b) = 0
		\end{eqnarray*}
		
		\item
		\begin{eqnarray*}
		D(a,b) = \sup \{ r \in \qp : \mathcal{A} \models R_{d(x, y) \geq r}(a, b) \} =_{(\ref{metric2})} \sup \{ r \in \qp : \mathcal{A} \models R_{d(x, y) \geq r}(b, a) \} = D(b, a)
		\end{eqnarray*}
		
		\item Let $c \in A$.  We want to show $D(a, c) \leq D(a, b) + D(b, c)$.  It is enough to show
		$$\forall r \in \qp (D(a, c) \geq r \implies D(a, b) + D(b, c) \geq r)$$
		Thus, let $r \in \qp$ and suppose $D(a, c) \geq r$.  Then $\sup \{ s \in \qp : \mathcal{A} \models R_{d(x, y) \geq s}(a, c) \} \geq r$.  By (\ref{metric3}), this means
		$$\sup \{ s \in \qp : \mathcal{A} \models \vee_{t \in \qp \cap [0,s]} R_{d(x, y) \geq t}(a, b) \wedge R_{d(x, y) \geq s-t}(b, c) \} \geq r$$
		Fix $n < \omega$.  There is some $s_n \in \qp$ such that $s_n \geq r - \frac{1}{n}$ and
		$$\mathcal{A} \models \vee_{t \in \qp \cap [0, s_n]} R_{d(x,y ) \geq t}(a, b) \wedge R_{d(x, y) \geq s_n -t}(b, c)$$
		Thus, there is some $t_n \in \qp$ such that $0 \leq t_n \leq s_n$ and
		$$\mathcal{A} \models R_{d(x, y) \geq t_n}(a, b) \wedge R_{d(x, y) \geq s_n - t_n}(b, c)$$
		By the definition of $D$, this means that $D(a, b) \geq t_n$ and $D(b, c) \geq s_n - t_n$; thus, $D(a, b) + D(b, c) \geq s_n$.  Since this is true for all $n < \omega$, we get that $D(a, b) + D(b, c) \geq r$ as desired.
	
	\end{enumerate}
	Thus, $D$ is a metric on $A$. \hfill $\dag_{Claim}$\\

Now we define functions and relations on $(A, D)$ such that they are uniformly continuous.  In particular,
\begin{enumerate}
	\item for $i < n_F$, set $f_i := \left(F_i^+\right)^\mathcal{A}$ with modulus 
	\begin{eqnarray*}
	\Delta_{f_i}(r) &=& \sup \{ s \in \qp : \mathcal{A} \vDash \forall \bx; \forall \by\\
	& &\left( \wedge_{i < n(F^+_i)} R_{d(z, z') \leq s}(x_i, y_i) \to R_{d(z, z') \leq r}(F^+_i(\bx), F^+_i(\by))\right) \}
	\end{eqnarray*}
	\item for $j < n_R$, set $r_j(\ba) := \sup \{ r \in \qp : \mathcal{A} \vDash R_{R_j(\bz) \leq r}[\ba] \}$ with modulus
	 \begin{eqnarray*}\Delta_{r_j}(r) &=& \sup \{ s \in \qp : \mathcal{A} \models \forall\bx \forall \by \\
	 & &\left(\wedge_{i < n(R_j)}R_{d(z, z') \leq s}(x_i, y_i) \to (R_{R_j(\bz)\monus R_j(\bz')\leq r}(\bx, \by) \wedge R_{R_j(\bz)\monus R_j(\bz')\leq r}(\by, \bx))\right) \}
	 \end{eqnarray*}
\end{enumerate}

Call this structure $\hat{\mathcal{A}}$.  We show that these functions and relations are uniformly continuous.  These moduli might not be the same moduli in the original signature $\tau$.  Instead, these are the optimal moduli, while the original language might have moduli that could be improved.

\begin{claim}
The functions $f_i$ and $r_j$ are uniformly continuous with their given moduli.
\end{claim}

{\bf Proof:} We show this for $f_i$; the proof for $r_j$ is similar.  Let $r \in \qp$ and let $\ba, \bb \in A$ such that $\max_{i < n} D(a_i, b_i) < \Delta_{f_i}(r)$.  Thus, for each $i < n$, $D(a_i, b_i) = \inf \{ s \in \qp : \mathcal{A} \models R_{d(x, y) \leq s}(a_i, b_i) \} < \Delta_{f_i}(r)$.  Since this is strict, there is some $s_i \in \qp$ such that $\mathcal{A} \models R_{d(x, y) \leq s_i}(a_i, b_i)$.  Note that (\ref{order2}) implies that the set $\Delta_{f_i}(r)$ is supremuming over is downward closed.  Thus, $s' = \max_{i < n} s_i$ is in it.  Thus, we can conclude
		$$\mathcal{A} \models R_{d(x, y) \leq r} [F^+_i(\ba), F^+_i(\bb)]$$
		This means that $D(f_i(\ba), f_i(\bb)) \leq r$, as desired. \hfill $\dag_{Claim}$\\
	
Now we have a prestructure.  Now we complete $\hat{\mathcal{A}}$ to $\overline{\mathcal{A}}$ in the standard way; see Munkries \cite{munkries} for a reference for the topological facts.  In particular, we define the continuous $\tau$-structure $\overline{\mathcal{A}}$ by
\begin{itemize}
	\item the universe $|\overline{\mathcal{A}}|$ is the completion of $(A, D)$;
	\item the metric $d^{\overline{\mathcal{A}}}$ is the extension of $D$ to $|\overline{\mathcal{A}}|$;
	\item for $i < n_F$, $F_i^{\overline{\mathcal{A}}}$ is the unique uniformly continuous extension of $f_i$ to $|\overline{\mathcal{A}}|$; and
	\item for $j < n_R$, $R_j^{\overline{\mathcal{A}}}$ is the unique uniformly continuous extension of $r_j$ to $|\overline{\mathcal{A}}|$.\\
\end{itemize}

We finish this part by showing that the desired structure is translated.

\begin{claim}
For all $\ba \in A$ and all formulas $\phi(\bx)$ built up from these functions and $D$, we have that
$$\phi^{\overline{\mathcal{A}}}(\ba) \square r \iff \mathcal{A} \vDash R_{\phi(\bz) \square r}[\ba]$$
\end{claim}

{\bf Proof:}  We proceed by induction on the construction of $\phi(\bx)$.  We assume that $\square$ is $\geq$ in our proofs, but the proofs for $\leq$ are the same.
\begin{itemize}
	
		\item If $\phi$ is atomic, then it falls into one of the following cases.
		\begin{itemize}

			\item Suppose $\phi(\bx) \equiv R_j(\term(\bx))$ for some term $\term$.  Then
			\begin{eqnarray*}
			R^{\overline{\mathcal{A}}}_j(\term(\ba)) \geq r &\iff& \inf \{ s \in \qp : \mathcal{A} \models R_{R_j(\bx) \geq s}[\term(\ba)] \} \geq r \\
			&\iff& \forall n < \omega, \exists s_n \in \qp \te{ so } s_n \geq r - \frac{1}{n} \te{ and } \mathcal{A} \models R_{R_j(\bx) \geq s_n}[\term(\ba)]\\
			&\iff_{(\ref{order2})}& \forall n < \omega, \mathcal{A} \models R_{R_j(\bx) \geq r - \frac{1}{n}}[\term(\ba)]\\
			&\iff_{(\ref{stdleft})}& \mathcal{A} \models R_{R_j(\bx) \geq r}[\term(\ba)]\\
			&\iff_{(\ref{termred})}& \mathcal{A} \models R_{R_j(\term(\by)) \geq r}[\ba]
			\end{eqnarray*}

			\item Suppose that $\phi(\bx, \by) \equiv d(\term_1(\bx), \term_2(\by))$ for terms $\term_1$ and $\term_2$.  The detail are essentially as above: $D^{\overline{\mathcal{A}}}(\term_1(\ba), \term_2(\bb))$ iff (by (\ref{order2}), the definition of sup, (\ref{stdright}) and (\ref{stdleft})) $\mathcal{A} \models R_{d(x, y) \geq r} [\term_1(\ba), \term_2(\bb)]$ iff (by (\ref{termred})) $\mathcal{A} \models R_{d(\term_1(\bx), \term_2(\by)) \geq r}(\ba, \bb)$.

		\end{itemize}
		
		\item For the inductive step, we deal with each connective (from our full set) in turn.  The induction steps for $x \mapsto 0$, $x \mapsto 1$, and $x \mapsto \frac{x}{2}$ are clear.
		\begin{itemize}
		
			\item Suppose $\phi \equiv \psi \monus \formula$.  Note if $r = 0$, then this is obvious.  So assume $r \neq 0$.  
			WLOG, we can assume $\psi^{\overline{\mathcal{A}}}(\ba) > \formula^{\overline{\mathcal{A}}}(\ba)$.

			\begin{itemize}
			
				\item First, suppose $\psi^{\overline{\mathcal{A}}}(\ba) - \formula^{\overline{\mathcal{A}}}(\ba) \geq r$.  Since $\psi^{\overline{\mathcal{A}}}(\ba) > \formula^{\overline{\mathcal{A}}}(\ba)$, there is some $s \in \qp$ such that $\psi^{\overline{\mathcal{A}}}(\ba) > s > \formula^{\overline{\mathcal{A}}}(\ba)$.  Then $\formula^{\overline{\mathcal{A}}}(\ba) \leq s$ and $\psi^{\overline{\mathcal{A}}}(\ba) \geq s + r$.  By induction, we have that
				$$\mathcal{A} \models R_{\formula(\bx) \leq s}[\ba] \wedge R_{\psi(\bx) \geq s + r}[\ba]$$
				Then, by (\ref{gmonus}), we have that $\mathcal{A} \models R_{\psi \monus \formula \geq r}[\ba]$ as desired.
				
				\item Now, suppose $\mathcal{A} \models R_{\psi \monus \formula \geq r}[\ba]$.  Again, (\ref{gmonus}) implies there is is some $s \in \qp$ such that
				$$\mathcal{A} \models R_{\formula \leq s}[\ba] \wedge R_{\psi \geq r+s}[\ba]$$
				By induction, we get $\formula^{\overline{\mathcal{A}}}(\ba) \leq s$ and $\psi^{\overline{\mathcal{A}}}(\ba) \geq r+s$.  Then
				$$\phi^{\overline{\mathcal{A}}}(\ba) = \psi^{\overline{\mathcal{A}}}(\ba) - \formula^{\overline{\mathcal{A}}}(\ba) \geq (r+s) - s = r$$
				as desired.
			
			\end{itemize}
			
			\item Suppose $\phi(\bx) \equiv \sup_x \psi(x, \bx)$.  We will consider both sides of the inequality since they're not symmetrically axiomatized (see (\ref{lsup}) and (\ref{gsup})). The case for $\inf$ is similar.
			\begin{itemize}
			
				\item Suppose that $\left(\sup_x \phi(x, \ba)\right)^{\overline{\mathcal{A}}} \geq r$.  Then for any $n < \omega$, there is some $a_n \in A$ such that $\phi^{\overline{\mathcal{A}}}(a_n, \ba) > r - \frac{1}{2n}$.  Since $\phi$ is uniformly continuous, there is some $\delta > 0$ such that, if $d(a_n, b) < \delta$, then 
				$$|\phi^{\overline{\mathcal{A}}}(a_n, \ba) - \phi^{\overline{\mathcal{A}}}(b, \ba)| < \frac{1}{2n}$$
				 Since $A$ is dense in $|\overline{\mathcal{A}}|$, there is some $a_n' \in A$ such that $d(a_n, a_n') < \delta$.  Thus, $\phi^{\overline{\mathcal{A}}}(a_n', \ba) > r - \frac{1}{n}$.  By induction, we have that
				$$\mathcal{A} \models \wedge_{n<\omega} \exists x R_{\phi(y, \by) \geq r - \frac{1}{n}}(x, \ba)$$
				Then (\ref{gsup}) says that $\mathcal{A} \models R_{\sup_y \phi(y, \by) \geq r}(\ba)$.
				
				\item Suppose that $\mathcal{A} \models R_{\sup_y \phi(y, \by) \geq r}[\ba]$.  Then, by (\ref{gsup}), 
				$$\mathcal{A} \models \wedge_{n<\omega} \exists x R_{\phi(y, \by) \geq r - \frac{1}{n}}(x, \ba)$$
				So, for each $n < \omega$, there is some $a_n \in A$ such that $\mathcal{A} \models R_{\phi(y, \by) \geq r - \frac{1}{n}}[a_n, \ba]$.  By induction, we have that $\phi^{\overline{\mathcal{A}}}(a_n, \ba) \geq r - \frac{1}{n}$.  Since this is true for each $n < \omega$, we get $\sup_y \phi^{\overline{\mathcal{A}}}(y, \ba) \geq r$.
				
				\item The other direction is easier and we can combine the two parts
				\begin{eqnarray*}
				\sup_x \phi^{\overline{\mathcal{A}}}(x, \ba) \leq r &\iff& \forall a \in |\overline{\mathcal{A}}|, \phi^{\overline{\mathcal{A}}}(a, \ba) \leq r \\
				&\iff& \forall a \in A, \phi^{\overline{\mathcal{A}}}(a, \ba) \leq r\\
				&\iff_{Induction}& \mathcal{A} \models \forall x R_{\phi \leq r}(x, \ba)\\
				&\iff_{(\ref{lsup})}& \mathcal{A} \models R_{\sup_x \phi(x, \bx)}[\ba]
				\end{eqnarray*}
			
			\end{itemize}
		
		\end{itemize}

	\end{itemize}
\hfill $\dag_{Claim}$\\

{\bf Essential inverses}\\

We now arrive at the final clause of Theorem \ref{goal}.  We need the following result to make it well-formed.

\begin{claim}
If $\mathcal{A} \vDash \dense_\tau$, then $A$ is nicely dense in $\overline{\mathcal{A}}$.
\end{claim}

{\bf Proof:} By construction, $A$ is dense in $\left|\overline{\mathcal{A}}\right|$.  For $i < n_F$, we have
$$F_i^{\overline{\mathcal{A}}} \rest A = \left(F^+_i\right)^{\mathcal{A}}$$
Thus, $A$ is closed under the functions of $\tau$ because $\mathcal{A}$ is a $\tau^+$-structure.\hfill$\dag_{Claim}$\\

\begin{prop} \label{ess-in-prop}
Given any continuous $\tau$-structure $M$ and nicely dense subset $A$, we have that $M \cong_A \overline{M_A}$ and, given any $\tau^+$ structure $\mathcal{A}$ that models $\dense_\tau$, we have that $\left(\overline{\mathcal{A}}\right)_A = A$.
\end{prop} \label{ess-in}

{\bf Proof:} First, let $M$ be a continuous $\tau$-structure and $A \subset M$ be nicely dense.  We define a map $f: M \to \overline{(M_A)}$ as follows: if $a \in A$, then $f(a) = a$.  For $a \in M-A$, fix some (any) sequence $\seq{a_n \in A : n < \omega}$ such that $\lim_{n \to \infty} a_n = a$ (this limit computed in $M$).  We know that $\seq{a_n : n < \omega}$ is Cauchy in $M$, so it's Cauchy in $\overline{(M_A)}$.  Then set $f(a) = \lim_{n \to \infty} a_n$, where that limit is computed in $\overline{(M_A)}$.  This is well-defined and a bijection because $A$ is dense in both sets.  That this is an $\tau$-isomorphism follows from applying the correspondence twice: for all $\ba \in A$ and $\phi(\bx) \in \cFml \tau$
$$\phi^M(\ba) \square r \iff M_A \vDash R_{\phi(\bx) \square r}[\ba] \iff \phi^{\overline{(M_A)}}(\ba) \square r$$
and the fact that the values of $\phi$ on $A$ determines its values on $M$ and $\overline{(M_A)}$.

Second, let $\mathcal{A}$ be a $\tau^+$ structure that models $\dense_\tau$.  Clearly, the universes are the same, i. e., $|(\overline{\mathcal{A}})_A| = A$.  For any relation $R_{\phi \square r}$ and $\ba \in A$, we have
$$\mathcal{A} \vDash R_{\phi \square r}[\ba] \iff \phi^{\overline{\mathcal{A}}}(\ba) \square r \iff \left(\overline{\mathcal{A}}\right)_A \vDash R_{\phi \square r}[\ba]$$
Given a function $F_i^+$ and $\ba, a \in A$, we have that
\begin{eqnarray*}
(F_i^+)^{\mathcal{A}}(\ba) = a &\iff& \mathcal{A} \vDash R_{d(F_i^+(\bx), x) \leq 0}[\ba, a] \\
&\iff& \left(\overline{\mathcal{A}}\right)_A \vDash R_{d(F_i^+(\bx), x) \leq 0}[\ba, a] \iff (F_i^+)^{\left(\overline{\mathcal{A}}\right)_A}(\ba) = a
\end{eqnarray*}
 \hfill $\dag_{Proposition\te{ }\ref{ess-in-prop},\te{ }Theorem\te{ }\ref{goal}}$\\

We can extend this correspondence to theories.  Suppose that $T$ is a continuous theory in $\tau$.  Following \cite[Definition 4.1]{fourguys}, theories are sets of closed $\tau$-conditions; that is, a set of ``$\phi=0$,'' where $\phi$ is a formula with no free variables.  The following is immediate from Theorem \ref{goal}.

\begin{cor}
If ``$\phi = 0$'' is a closed $\tau$-condition, then
$$\phi^M=0 \iff M_A \vDash R_{\phi \leq 0}$$
\end{cor}

With our fixed theory $T$, set $T^*$ to be $\dense_\tau \cup \{ R_{\phi\leq 0} : ``\phi=0'' \in T\}$.  Then our representation of continuous $\tau$-structures as discrete $\tau^+$-structures modeling $\dense_\tau$ can be extended to a representation of continuous models of $T$ and discrete models of $T^*$.

\section{Elementary Substructure} \label{precsec}

We now discuss translating the notion of elementary substructure between our two contexts.  Depending on the generality needed, this is either easy or difficult.

For the easy case, we have the following.

\begin{theorem}\label{4.1}
Let $M, N$ be continuous $\tau$-structures.  The following are equivalent:
\begin{enumerate}
	\item $M \prec_\tau N$.
	\item For all nicely dense $A \subset M$ and $B \subset N$ such that $A \subset B$, we have $M_A \subset_{\tau^+} N_B$.
	\item There exist nicely dense $A \subset M$ and $B \subset N$ such that $A \subset B$ and $M_A \subset_{\tau^+} N_B$.
\end{enumerate}Then $M \prec_\tau N$ iff, for every nicely dense $A \subset M$ and $B \subset N$ such that $A \subset B$, we have that $M_A \subset_{\tau^+} N_B$.
\end{theorem}

Note that the relation between $M_A$ and $N_B$ is just substructure.  So even though they are models of infinitary theories, their relation just concerns atomic formulas.  This is because we have built the quantifiers of $\tau$ into the relations of $\tau^+$.

{\bf Proof:} $(1) \to (2)$:  Let $A \subset M$ and $B \subset N$ be nicely dense so $A \subset B$.  We want to show that $M_A \subset_{\tau^+} N_B$.
\begin{itemize}
	\item Let $F^+ \in \tau^+$ and $\ba \in A$.  Then, by definition of the structures,
	$$(F^+)^{M_A}(\ba) = F^M(\ba) = F^N(\ba) = (F^+)^{N_B}(\ba)$$
	\item Let $R_{\phi \square r} (\bx) \in \tau^+$ and $\ba \in A$.
	\begin{align*}
	M_M \vDash R_{\phi(\bx) \square r}[\ba] &\iff \phi^M(\ba) \square r & \te{Theorem }\ref{goal}\\
	&\iff \phi^N(\ba) \square r & M\prec_\tau N\\
	&\iff N_N \vDash R_{\phi(\bx) \square r}[\ba] & \te{Theorem }\ref{goal}
	\end{align*}	
\end{itemize}

$(2) \to (3)$: Immediate.

$(3) \to (1)$: We want to show that, for all $\ba \in M$ and $\phi(\bx) \in \cFml \tau$, we have $\phi^M(\ba) = \phi^N(\ba)$.  By the continuity of formulas, it is enough to show this for a dense subset of $M$, namely $A$.  So let $\ba \in M$.  From the theorems proved last section, we have, for each $r \in \qp$,
\begin{align*}
\phi^M(\ba) \square r &\iff M_A \vDash R_{\phi(\bx) \square r}[\ba] & \te{Theorem }\ref{goal}\\
&\iff N_B \vDash R_{\phi(\bx) \square r}[\ba] & M_A \subset_{\tau^+} N_B\\
&\iff \phi^N(\ba) \square r & \te{Theorem }\ref{goal}
\end{align*}
Thus $\phi^M(\ba) = \phi^N(\ba)$ and $M\prec_\tau N$ as desired. \hfill\dag\\

This gives us the following corollary.

\begin{cor}\label{foraec}
Given $\mathcal{A}, \mathcal{B} \vDash \dense_\tau$ with $A \subset B$, $\mathcal{A} \subset_{\tau^+} \mathcal{B}$ iff $\overline{\mathcal{A}} \prec_{\tau} \overline{\mathcal{B}}$.
\end{cor}

Corollary \ref{foraec} characterizes when $\tau^+$-structures $\mathcal{A}, \mathcal{B}$ give rise to $\tau$-structures that are related by elementary substructure \emph{if} $A \subset B$.  However, we would also like to know when this happens without the assumption that $A$ is a subset of $B$.  A natural example of when this happens is that $\mathbb{Q} \cap[0,1]$ and $\mathbb{Q} +\sqrt{2}$ are disjoint, but their completions fit together nicely.  Theorem \ref{4.1} characterizes when $\overline{\mathcal{A}} \prec_\tau \overline{\mathcal{B}}$ by the existence of nicely dense $C \subset \left|\overline{\mathcal{A}}\right|$ and $D\subset \left|\overline{\mathcal{B}}\right|$ such that $C \subset D$ and $\left( \overline{\mathcal{A}}\right)_C \subset_{\tau^+} \left(\overline{\mathcal{B}}\right)_D$.

This is not ideal because it requires quantifying over subsets of the completion of $\tau^+$-structure.  We would prefer a characterization more closely tied to just the $\tau^+$ structure on $\mathcal{A}$ and $\mathcal{B}$.  To this end, we define inessential extensions of models of $\dense_\tau$.

\begin{defin}
Given $\mathcal{A} \subset_{\tau^+} \mathcal{B}$, we say that $\mathcal{B}$ is an \emph{inessential extension} of $\mathcal{A}$ iff for every $b \in B$ and $n < \omega$, there is some $a \in A$ such that $\mathcal{B} \vDash R_{d(x, y) < \frac{1}{n}}[b, a]$.
\end{defin}

Briefly, we have that $\mathcal{B}$ is an inessential extension of $\mathcal{A}$ if $A \subset B$ and $\overline{\mathcal{A}} = \overline{\mathcal{B}}$.  This gives us our desired criterion.

\begin{theorem}
Let $\mathcal{A}, \mathcal{B} \vDash \dense_\tau$.  The following are equivalent:
\begin{enumerate}
	\item $\overline{\mathcal{A}} \prec_\tau \overline{\mathcal{B}}$.
	
	\item There is a $\tau^+$-structure $\mathcal{C}\vDash \dense_\tau$ such that $\mathcal{A}, \mathcal{B} \subset_{\tau^+} \mathcal{C}$ and $\mathcal{C}$ is an inessential extension of $\mathcal{B}$.

	\item There is a $\tau^+$-structure $\mathcal{A} \cup \mathcal{B} \vDash \dense_\tau$ such that $\mathcal{A}, \mathcal{B} \subset_{\tau^+} \mathcal{A} \cup \mathcal{B}$; the universe of $\mathcal{A} \cup \mathcal{B}$ is the closure of $A \cup B$ under functions; and such that for every $a \in |\mathcal{A} \cup \mathcal{B}|$ and $n < \omega$, there is some $b \in B$ such that $\mathcal{A} \cup \mathcal{B} \vDash R_{d(x, y) < \frac{1}{n}}[b, a]$.
\end{enumerate}
\end{theorem}

{\bf Proof:} $(1) \to (2)$:  Let $C$ be the closure of $A \cup B$ under the functions of $\overline{\mathcal{B}}$.  Then $C$ is a nicely dense subset of $\overline{\mathcal{B}}$.  Take $\mathcal{C} = \left(\overline{\mathcal{B}}\right)_{C}$.

$(2) \to (3)$: Set $\mathcal{A} \cup \mathcal{B}$ to be the $\tau^+$-substructure of $\mathcal{C}$ whose universe is generated from $A \cup B$.

$(3) \to (1)$: We have $\overline{\mathcal{A}}, \overline{\mathcal{B}} \prec \overline{\mathcal{A} \cup \mathcal{B}}$ from the first condition and $\overline{\mathcal{B}} = \overline{\mathcal{A}\cup \mathcal{B}}$ from the second.
\hfill\dag\\
\section{$\dense_\tau$ as an Abstract Elementary Class}

In this section, we view the discrete side of things as an Abstract Elementary Class; see Baldwin \cite{baldwinbook} or Grossberg \cite{ramibook}.

\begin{theorem} \label{Tdense-aec}
Let $T$ be a complete, continuous first order $\tau$-theory.  Then let $\tau^+$ and $\dense_\tau$ be from Theorem \ref{goal}.  Set $K = (\Mod (T^*), \subset_{\tau^+})$.  Then
\begin{enumerate}
	\item $K$ is an AEC;
	\item $K$ has amalgamation over sets, joint embedding, and no maximal models; and
	\item Galois types in $K$ correspond to continuous, syntactic types in $T$.
\end{enumerate} 
\end{theorem}

Note that if $T$ were not complete, then amalgamation over sets would not hold.  However, the other properties will continue to hold, including the correspondence between Galois types and syntactic types.

{\bf Proof:}  $T^*$ is a $\mathbb{L}_{\omega_1, \omega}(\tau^+)$ theory, so all of the axioms of AECs hold except perhaps the chain axioms, which might fail because strong substructure is not elementarity according to the fragment.  

Let $\seq{\mathcal{A}_i:i < \alpha}$ be a $\subset_{\tau^+}$-increasing chain of models of $T^*$.  Then, by Theorem \ref{goal} and Corollary \ref{foraec}, the sequence $\seq{\overline{\mathcal{A}_i} : i < \alpha}$ is $\prec_\tau$-increasing chain that each model $T$.  Then by the chain axiom for continuous logic, there is $M = \overline{\cup_{i < \alpha} (\overline{\mathcal{A}_i})}$ that models $T$.  Additionally, $A:=\cup_{i < \alpha} A_i$ is nicely dense in $M$.  Here, $A$ is closed under functions because they are finitary and $\alpha$ is limit.  Thus, $M_A = \cup_{i<\alpha} \mathcal{A}_i$ is as desired.  Additionally, suppose $\mathcal{A}_i \subset_{\tau^+} \mathcal{B}$ for all $i < \alpha$.  Then we have $M \prec_\tau \overline{\mathcal{B}}$, so $M_{A} \subset_{\tau^+} \mathcal{B}$.

The rest of these properties all follow from the corresponding properties of continuous first-order logic.  For instance, considering amalgamation, suppose $\mathcal{A} \subset_{\tau^+} \mathcal{B}, \mathcal{C}$.  Then we have $\overline{\mathcal{A}} \prec_\tau \overline{\mathcal{B}}, \overline{\mathcal{C}}$.  By amalgamation for continuous first-order logic, there is some $N \succ_\tau \overline{\mathcal{B}}$ and elementary $f:\overline{\mathcal{C}} \to_{\overline{\mathcal{A}}} N$.  Let $D \subset N$ be nicely dense that contains $B \cup C$.  Then we have $\mathcal{B} \subset_{\tau^+} N_D$ and $f\rest C:\mathcal{C} \to_{\mathcal{A}} N_D$; this is an amalgamation of the original system.

Finally, we wish to show that Galois types are syntactic types and vice versa.  Note that there are monster models in each class.  Further more, we may assume that, if $\sea$ is the monster model of $T$, then there is some nicely dense $U \subset \sea$ such that the monster model of $K$ is $M_U$.  In fact, we could take $U = |\sea|$.  Let continuous $M \vDash T$ and $A \subset M$ be nicely dense.  If we have tuples $\ba$ and $\bb$, then
\begin{eqnarray*}
gtp_K(\ba/M_A) = gtp_K(\bb/M_A) &\iff& \exists f \in Aut_{M_A} M_U. f(\ba) = \bb\\
&\iff& \exists f \in Aut_{\overline{M_A}} \sea. f(\ba) = \bb\\
&\iff& tp_T(\ba/M) = tp_T(\bb/M)
\end{eqnarray*}
\hfill \dag\\

We pause here only briefly to point out a strange occurence: first-order continuous logic is compact (see \cite[Theorem 5.8]{fourguys}), but $\Ll_{\omega_1, \omega}$ is incompact.  Yet, we have seen that continuous logic can be embedded into $\Ll_{\omega_1, \omega}$.  The solution to this incongruity is that the compactness of continuous logic comes from a different ultraproduct than the model-theoretic one, namely the metric space ultraproduct.  It is possible to formulate this ultraproduct within the context of models of $\dense_\tau$ by removing functions that would lead to nonstandard distances and changing the notion of equivalence.  This is done for the case of Banach spaces in \cite[Section 4.1]{gam-up}.

We conclude with an example relevant to the next section.  I thank Ilijas Farah for discussions about this example.

\begin{example}
The model theory of probability spaces $(X, B, \mu)$ is developed in \cite[Section 16]{fourguys} by defining a metric structure where the universe consists of the elements of $B$ modulo $\mu$-equivalence; adding 0, 1, complement, meet, and join to the language; and setting the distance $d(A, B) = \mu (A \triangle B)$. An axiomatization $PrA$ is given in \cite[Theorem 16.1]{fourguys}.

Then $PrA^*$ consists of the normal axioms of $\dense_\tau$ along with the appropriately relativized axioms of $PrA$: each ``universal axiom'' $\sup_\bx \phi(\bx)$ becomes $\forall \bx R_{\phi\leq 0}(\bx)$.  This means that a model $M$ of $PrA^*$ will have a naturally induced structure of a boolean algebra on it, along with a distance function $d$ that gives rise to a probability measure via $\mu(x) = d(x, 0)$.
\end{example}

Probability spaces can be seen as a metric version of boolean algebras, and any model of $PrA^*$ is a boolean algebra with some extra structure.  However, this correspondence does not go backwards: there are many boolean algebras that cannot be expanded to a model of $PrA^*$.  In particular, any boolean algebra $\mathcal{B}$ with an increasing sequence $\seq{b_\alpha : \alpha < \omega_1}$ cannot be a reduct of a model of $PrA^*$ because $\seq{\mu(b_\alpha) : \alpha<\omega_1}$ would be $\aleph_1$-increasing sequence of reals.

$PrA$ is stable (see \cite[Proposition 16.9]{fourguys} or \cite[Proposition 4.6]{mtoI}) and, thus, doesn't have the (properly defined) order property.  However, there \emph{is} an order extractable from the syntax.  Set $x \prec y$ iff $\mu(x \cap y) \monus \mu(x) = 0$.  Looking in the the probability space of Lebesgue measurable subsets of $[0,1]$, we can find sequences $\seq{X_i \mid i< \alpha}$ that are increasing  according to this measure for any $\alpha < \omega_1$.  However, this does not contradict the above because the proper definition of the order property \cite[Definition 2.3]{mtoI} requires a uniform separation of the values of the formula (rather than only detecting when a formula takes the value 0).

This represents a common theme of studying the order property in general AECs.  When defining the order property in an AEC (see, for instance, \cite[Section 4]{sh394}), an extra parameter of \emph{how long} the order is must be included as there is no compactness to make it arbitrarily long.  Then the order property with no parameters means that there are arbitrarily long orders or, equivalently, that there is an order of length $\beth_{(2^{LS(K)})^+}$.  Thus, $PrA^*$ has the $\alpha$-order property for every $\alpha < \omega_1$, does not have the $\aleph_1$-order property, and is stable (see Corollary \ref{stab-trans-cor} below).  The compactness of continuous logic from the metric ultraproduct offers a better criteria in terms uniform order property from \cite{mtoI}: $T^*$ has the order property (as an AEC) iff $T$ has arbitrary long finite orders (as a metric theory).

\section{Types} \label{type-sec}

Fix a continuous $\tau$-theory $T$, and let $T^*$ be its discrete equivalent (described at the end of Section \ref{main-sect}).

In Theorem \ref{Tdense-aec}, we showed that Galois types on the discrete side correspond to the normal definition of syntactic on the continuous side.  We can characterize Galois types more precisely by looking at the quantifier-free syntactic type in $\tau^+$.

\begin{prop}
In $(\Mod (T^*), \subset_{\tau^+})$, two elements have the same Galois type iff they have the same quantifier-free $\tau^+$-type.
\end{prop}

{\bf Proof:} In all AECs, two tuples having the same Galois type implies having the same quantifier-free type.  Suppose that $\ba, \bb \in \mathcal{B} \vDash T^*$ have the same quantifier-free type over $A \subset B$. Then, by Theorem \ref{goal}.(2), $\ba, \bb \in \overline{\mathcal{B}}$ have the same continuous $\tau$-type over $C$.  By Theorem \ref{Tdense-aec}.(3), this means that $\ba$ and $\bb$ have the same Galois type in $( \Mod(T^*), \subset_{\tau^+})$.\hfill \dag\\

Of course, not every finitely consistent quantifier-free $\tau^+$-type is realized in a model.  For instance, the partial type $\{ \neg R_{d(x',y') \leq 0}(x,y), R_{d(x',y') \leq \frac{1}{n}}(x, y) \mid n < \omega\}$ is finitely satisfiable, but describes the type of two elements that are an infinitesimal, positive distance from each other.  Clearly, no element in a model of $\dense_\tau$ can satisfy this.

The correspondence allows us to connect the stability of the discretization with the stability of the continuous class because we've established a bijection between the types.  An important distinction in the study of stability in continuous first-order logic is \emph{discrete stability} (counting the cardinality of the type space \cite[Definition 14.1]{fourguys}) and \emph{stability} (counting the density character of the type space under a natural metric \cite[Definition 14.4]{fourguys}).

\begin{cor} \label{stab-trans-cor}
$T$ is $\lambda$-discrete stable iff $(\Mod (T^*), \subset_{\tau^+})$ is $\lambda$-Galois stable.
\end{cor}

{\bf Proof:} Theorem \ref{Tdense-aec}.(3) has established a bijection between the relevant sets.\hfill\dag\\

The density character notion of stability turns out to be the better measure of a class's behavior (although they are equivalent globally by \cite[Theorem 14.6]{fourguys}).  Additionally, in the discrete context, we could develop a notion of distance between types paralleling the notion on the continuous side.  However, it doesn't seem as though there's a natural way to recover the density character of the type space in this context.

We now turn to saturation.  While stability is a property of theories, saturation is a property of models.  Thus, it requires a finer analysis to determine the transfer between continuous and discrete notions.  In particular, it is clear that we cannot expect a characterization along the lines of ``the continuous structure $M$ is saturated iff the discrete structure $M_A$ is saturated for every nicely dense subset $A$."  The failure comes from the fact that nicely dense subsets can miss many elements.  Thus, we introduce the notion of a sequence type below.  This notion allows us to give a characterization of the desired form in terms of saturation for sequence types (see Theorem \ref{Tdense-sat}).  Note this notion appears elsewhere, e.g., in Farah and Magidor \cite{farah-magidor} as the type $p_\omega$, although they formalize it as a type in infinitely many free variables, rather than infinitely many types in 2 free variables.

Recall that, in continuous first-order logic, formulas are also uniformly continuous functions with a modulus computable from the language.  Given $\phi \in \cFml \tau$, write $\Delta_\phi$ for this modulus.

\begin{defin}\label{sequencetypedef}\
\begin{itemize}

	\item We say that $\seq{r_n : n < \omega}$ is a \emph{sequence $\ell$-type over $A$} iff $r_0(\bx)$ is an $\ell$-type over $A$ and $r_{n+1}(\bx, \by)$ is a $2 \ell$-type over $A$ such that there is some index set $I$, (possibly repeating) formulas $\{\phi_i : i \in I\}$; and (possibly repeating) Cauchy sequences $\{ \seq{\bb_n^i \in B'}_{n<\omega} : i \in I\}$ so $d(\bb_n^i, \bb_{n+1}^i) \leq \frac{1}{2^n}$ such that
	\begin{itemize}
		\item $r_0(\bx) = \{ R_{\phi_i(\bz, \bz') \leq \Delta_{\phi_i}(2)}(\bx, \bb^i_0): i \in I\}$; and
		\item $r_{n+1}(\bx, \by) = \{ R_{\phi_i(\bz, \bz') \leq \Delta_{\phi_i}(\frac{1}{2^n})}(\bx, \bb^i_{n+1}): i \in I\} \cup \{ R_{d(z, z') \leq \frac{1}{2^n}}(x_k, y_k) : k < \ell \}$
	\end{itemize}

	\item A \emph{realization} of a sequence type $\seq{r_n : n <\omega}$ is $\seq{\ba_n : n < \omega}$ such that
	\begin{itemize}
		\item $\ba_0$ realizes $r_0$; and
		\item $\ba_{n+1}\ba_n$ realizes $r_{n+1}$.
	\end{itemize}
	
\end{itemize}
\end{defin}

Note that the use of $\frac{1}{2^n}$ is not necessary; this could be replaced by any summable sequence for an equivalent definition (also replacing $\frac{1}{2^{n-1}}$ by the trailing sums).  However, we fix $\frac{1}{2^n}$ for computational ease.  The fundamental connection between continuous types and sequence types is the following.

\begin{theorem}\label{limittypes}
Let $A \subset M$ be nicely dense.
\begin{enumerate}

	\item If $B \subset M$ and $r(\bx)$ is a partial $\ell$-type over $B$, then for any $B' \subset A$ such that $\overline{B'} \supset B$, there is a sequence $\ell$ type $\seq{r_n:n<\omega}$ over $B'$ such that 
	\begin{center}
	$M$ realizes $r$ iff $M_A$ realizes $\seq{r_n : n < \omega}$
	\end{center}

	\item If $B' \subset A$ and $\seq{r_n :n < \omega}$ is a partial sequence $\ell$-type over $B'$, then there is a unique $\ell$-type $r$ over $B'$ such that
	\begin{center}
	$M$ realizes $r$ iff $M_A$ realizes $\seq{r_n : n < \omega}$
	\end{center}
	
\end{enumerate}
\end{theorem}

In (2), $r$ extends to a unique type over $\overline{B'}$.  We sometimes write $\lim_{n\to\infty} r_n$ for $r$.  In each case, we have that $\seq{\ba_n \in M_A : n < \omega}$ realizes $\seq{r_n : n < \omega}$ implies $\lim_{n\to\infty} \ba_n$ realizes $\lim_{n\to\infty} r_n$.

{\bf Proof:}\begin{enumerate}

	\item For $n < \omega$ and $b \in B$, set $B'_n(b) = \{ b' \in B' : d^M(b', b) < \frac{1}{2^n} \}$.  To make the cardinality work out nicer, fix a choice function $G$, i. e., $G(B'_n(b)) \in B'_n(b)$.  Then $B'_n(\bb)$ and $G(B'_n(\bb))$ have the obvious meanings.  Define
\begin{eqnarray*}
r_n^+(\bx) &:=& \{ R_{\phi(\bz; \by) \leq \Delta_\phi(\frac{1}{2^{n-1}})}(\bx; G(B'_n(\bb))) : \te{``}\phi(\bx; \bb) = 0\te{''} \in r \}\\
r_0(\bx) &:=& r^+_0(\bx)\\
r_{n+1}(\bx, \by) &:=& r^+_{n+1}(\bx) \cup \{ R_{d(z, z') \leq \frac{1}{2^n}}(x_i, y_i) : i < \ell(\bx) \}
\end{eqnarray*}
Then $\seq{r_n:n<\omega}$ is a sequence type over $B'$; we can see this by taking $r$ as the index set, $\phi_i = \phi$, and $\bb^i = G(B'_n(\bb))$ for $i = \te{``}\phi(\bx; \bb) = 0\te{''} \in r$.  

To show it has the desired property, first suppose that $\seq{\ba_n:n < \omega}$ from $M_A$ realizes $\seq{r_n:n<\omega}$.  We know that $\seq{\ba_n : n < \omega}$ is a Cauchy sequence; in particular, for $m > n$, $$d^M(\ba_n, \ba_m) \leq \sum_{i = n}^{m} \frac{1}{2^i} = \frac{2^{m+1-n}-1}{2^{m+1}} $$
Since $\overline{M_A} \cong M$ is complete, there is $\ba \in M$ such that $\lim_{n\to\infty}\ba_n = \ba$.  We claim that $\ba \vDash r$.  Let ``$\phi(\bx; \bb) = 0\te{''}\in r$.  Then
$$d^M\left(\ba\bb; \ba_n G(B'_n(\bb))\right) \leq \max \left\{ \sum_{i=n}^\infty \frac{1}{2^i}, \frac{1}{2^n} \right\} = \frac{1}{2^{n-1}}$$
Thus, 
$$|\phi^{M}(\ba;\bb) - \phi^M(\ba_n; G(B'_n(\bb)))| \leq \Delta_\phi\left(\frac{1}{2^{n-1}}\right)$$
Letting $n \to \infty$, we have that 
$$\phi^M(\ba;\bb) = \lim_{n \to \infty} \phi^M(\ba_n; G(B'_n(\bb))) \leq \lim_{n \to \infty} \Delta_\phi\left(\frac{1}{2^{n-1}}\right)= 0$$
as desired.

Second, suppose that $\ba \in M$ realizes $r$.  Since $A$ is dense, we can find $\seq{\ba_n \in A : n < \omega}$ such that $d^M(\ba_n, \ba_{n+1}) < \frac{1}{2^n}$ and $\ba_n \to \ba$.  We want to show that $r_n^+(\ba_n)$ holds.  Let ``$\phi(\bx, \bb)\te{''}\in r$.  We know that
$$d^M(\ba_n G(B'_n(\bb)), \ba\bb) \leq \frac{1}{2^{n-1}}$$
Thus,
$$\phi^{M_A}(\ba_n, G(B'_n(\bb))) = \left| \phi^M(\ba, \bb) - \phi^M(\ba_n, G(B'_n(\bb))) \right| < \Delta_\phi\left(\frac{1}{2^{n-1}}\right)$$
as desired.\hfill \dag\\

	\item  Let $\seq{r_n : n < \omega}$ be a partial sequence $\ell$-type given by $I$, $\{\phi_i :i \in I\}$, and $\{\seq{\bb_n^i}_{n<\omega}:i\in I\}$.  Then set
	$$r(\bx):=\left\{ \phi_i(\bx,  \bb_n^i) \leq 2 \Delta_{\phi_i}\left(\frac{1}{2^{n-1}}\right) : i \in I\right\}$$
	First, suppose that $\seq{\ba_n \in M_A:n < \omega}$ realizes $\seq{r_n: n < \omega}$.  This is a Cauchy sequence, so we can find its limit $\ba \in M$.  We have that $d^M(\ba\bb^i_n, \ba_n\bb_n) \leq \frac{1}{2^{n-1}}$.  Thus, by the definition of moduli,
	\begin{eqnarray*}
	\left| \phi^M_i (\ba, \bb_n^i) - \phi^M_i (\ba_n, \bb_n) \right| \leq \Delta_{\phi_i}\left(\frac{1}{2^{n-1}}\right)\\
	\phi^M_i(\ba, \bb_n^i) \leq \Delta_{\phi_i}\left(\frac{1}{2^{n-1}}\right) + \phi^M_i (\ba_n, \bb_n) \leq 2\Delta_{\phi_i}\left(\frac{1}{2^{n-1}}\right)
	\end{eqnarray*}
So $\ba \vDash r$.

Second, suppose that $\ba \in M$ realizes $r$.  By continuity, this implies $\phi^M_i(\ba, \lim_{n\to\infty}\bb^i_n) = 0$.  Then, by denseness, we can find a Cauchy sequence $\seq{\ba_n \in A:n<\omega}$ such that $d(\ba_{n+1}, \ba_n) \leq \frac{1}{2^n}$.  Then $d(\ba \lim_{n \to \infty} \bb_n^i, \ba_n \bb_n^i) \leq \frac{1}{2^{n-1}}$.  Then we can conclude

$$|\phi^M_i(\ba, \lim_{n \to \infty} \bb_n^i) - \phi^M_i(\ba_n, \bb_n^i)| \leq \Delta_{\phi_i}(\frac{1}{2^{n-1}})$$
$$\phi_i^{M}(\ba_n, \bb_n^i) \leq \Delta_{\phi_i}(\frac{1}{2^{n-1}})$$
$$M_A \vDash R_{\phi_i(\bx, \by) \leq \Delta_{\phi_i}\left(\frac{1}{2^{n-1}}\right)} [\ba_n, \bb_n^i]$$

So $\seq{\ba_n:n<\omega}$ realizes $\seq{r_n:n<\omega}$.

\end{enumerate}\hfill \dag\\

We now connect type-theoretic concepts in continuous logic (e.g. saturation and stability) with concepts in our discrete analogue.

Recall (see \cite[Definition 7.5]{fourguys}) that a continuous structure $M$ is $\kappa$-saturated iff, for any $A \subset M$ of size $< \kappa$ and any continuous type $r(\bx)$ over $A$, if every finite subset of $r(\bx)$ is satisfiable in $M$, then so is $r(\bx)$.

\begin{defin} \
\begin{itemize}
	\item If $\seq{r_n:n<\omega}$ is a sequence type defined by an index set $I$ and $I_0 \subset I$, then $\seq{r_n : n < \omega}^{I_0}$ is the sequence type defined by $I_0$.
	\item We say that $M_A \vDash \dense_\tau$ is $\kappa$-saturated for sequence types iff, for all $B' \subset A$ of size $<\kappa$ and sequence type $\seq{r_n:n<\omega}$ over $B'$ that is defined by $I$, if $\seq{r_n:n<\omega}^{I_0}$ is realized in $M_A$ for all finite $I_0 \subset I$, then $\seq{r_n : n < \omega}$ is realized in $M_A$.
\end{itemize}
\end{defin}

We can characterize the saturation of $M$ by the saturation for sequence types of an $M_A$.

\begin{theorem} \label{Tdense-sat} Let $M \vDash T$ and $A \subset M$ be nicely dense.  Then $M$ is $\kappa$-saturated iff $M_A$ is $\kappa$-saturated for sequence types.
\end{theorem}

{\bf Proof:} First, let $M$ be $\kappa$-saturated and $A \subset M$ be nicely dense.  Let $B' \subset M_A$ of size $\lambda$ and let $\seq{r_n : n < \omega}$ be a sequence type over $B'$ that is finitely satisfiable in $M_A$.  Set $r=\lim_{n\to\infty}r_n$ from Theorem \ref{limittypes}; this is a type over $B'$.  We claim that $r$ is finitely satisfiable in $M$.  Any finite subset of $r_0$ of $r$ corresponds to a finite $I_0 \subset I$.  Then, by Theorem \ref{limittypes}, $r_0$ is realized in $M$ iff $\seq{r_n:n<\omega}^{I_0}$ is realized in $M_A$.  Then, since each $\seq{r_n:n<\omega}^{I_0}$ is realized in $M_A$ by assumption, we have that $r$ is finitely satisfiable in $M$.  By the $\kappa$-saturation of $M$, $r$ is realized in $M$.  By Theorem \ref{limittypes}, $\seq{r_n:n<\omega}$ is realized in $M_A$.  So $M_A$ is $\kappa$-saturated for sequence types.
	
Let $M_A$ be $\kappa$-saturated for sequence types.  Let $B \subset M$ of size $<\kappa$ and $r$ be a type over $B$ that is finitely satisfied in $B$.  Find $B' \subset A$ such that $\overline{B'} \supset B$; this can be done with $|B'|\leq |B|+\aleph_0 <\kappa$.  Then form the sequence type $\seq{r_n:n<\omega}$ over $B'$ that converges to $r_n$ as in Theorem \ref{limittypes}.  As before, since $r$ is finitely satisfiable in $M$, so is $\seq{r_n:n<\omega}$ in $M_A$.  So $\seq{r_n:n<\omega}$ is realized in $M_A$ by saturation.  Thus, $r$ is realized in $M$. \hfill \dag\\

\section{Metric Abstract Elementary Classes}

In this section, we extend the above representation to Metric Abstract Elementary Classes.  Recall from Hirvonen and Hyttinen \cite{maecs} or that a Metric Abstract Elementary Class (MAEC) is a class of continuous $\tau(K)$-structures $K$ and a strong substructure relation $\prec_K$ satisfying the following axioms:
\begin{enumerate}

    \item $\prec_K$ is a partial order on $K$;

    \item for every $M, N \in K$, if $M \prec_K N$, then $M \subseteq_{\tau} N$;

    \item $(K, \prec_K)$ respects $\tau(K)$ isomorphisms, if $f: N \to N'$ is an $\tau(K)$ isomorphism and $N \in K$, then $N' \in K$ and if we also have $M \in K$ with $M \prec_K N$, then $f(M) \in K$ and $f(M) \prec_K N'$;

    \item \emph{(Coherence)} if $M_0, M_1, M_2 \in K$ with $M_0 \prec_K M_2$; $M_1 \prec_K M_2$; and $M_0 \subseteq M_1$, then $M_0 \prec M_1$;

    \item \emph{(Tarski-Vaught chain axioms)} suppose $\seq{M_i \in K : i < \alpha}$ is a $\prec_K$-increasing continuous chain, then

        \begin{enumerate}

            \item $\overline{\cup_{i < \alpha} M_i} \in K$ and, for all $i < \alpha$, we have $M_i \prec_K \overline{\cup_{i < \alpha} M_i}$; and

            \item if there is some $N \in K$ such that, for all $i < \alpha$, we have $M_i \prec_K N$, then we also have $\cup_{i < \alpha} M_i \prec_K N$; and

        \end{enumerate}

    \item \emph{(L\"{o}wenheim-Skolem number)} There is an infinite cardinal $\lambda \geq |\tau(K)|$ such that for any $M \in K$ and $A \subset M$, there is some $N \prec_K M$ such that $A \subset |N|$ and $dc(N) \leq |A| + \lambda$.  We denote the minimum such cardinal by $LS(K)$.

\end{enumerate}

These axioms were first given in Hirvonen and Hyttinen \cite{maecs}.

A key difference is that the functions and relations $\tau$ are no longer required to be \emph{uniformly} continuous, but just continuous.  
This is due to the lack of compactness in the MAEC context.  This initially seems problematic because functions must be uniformly continuous on a set to be guaranteed an extension to its closure.  However, we get around this by simply defining $K_{dense}$ to be all the structures that happen to complete to a member of $K$, then use the MAEC axioms to show that $K_{dense}$ satisfies the AEC axioms.\\

For this reason, when we refer to continuous languages, structures, etc. in this section, we will not mean that they are uniformly continuous.

\begin{theorem} \label{Kdense}
Let $\tau$ be a continuous language.  Then there is a discrete language $\tau^+$ such that, for every MAEC $K$ with $\tau(K)= \tau$, there is
\begin{enumerate}
	\item an AEC $K_{dense}$ with $\tau(K_{dense}) = \tau^+$ and $LS(K_{dense}) = LS(K)$;
	\item a map from $M \in K$ and nicely dense subsets $A$ of $M$ to $M_A \in K_{dense}$; and
	\item a map from $\mathcal{A} \in K_{dense}$ to $\overline{\mathcal{A}} \in K$
\end{enumerate}
with the properties that
\begin{enumerate}
	\item $M_A$ has universe $A$ and for each $\ba \in A$, $r \in \qp$, and $\square \in \{\leq, \geq\}$, we have that 
\begin{eqnarray*}
M_A \vDash R_{R_j(\bz) \square r}[\ba] \iff R^M_j(\ba)\square r
\end{eqnarray*}
\item $\overline{\mathcal{A}}$ has universe that is the completion of $A$ with respect to the derived metric and for each $\ba \in A$, $r \in \qp$, and $\square \in \{\leq, \geq\}$, we have that
$$R_j^{\overline{\mathcal{A}}}(\ba) \square r \iff \mathcal{A} \vDash R_{R_j(\bz) \square r}[\ba]$$
\item The maps above are essentially inverses, in the since of Theorem \ref{goal}.3
\item \begin{itemize}
	\item Given $M_\ell \in K$ and $A_\ell$ nicely dense in $M_\ell$ for $\ell = 0,1$, if $f:M_0 \to M_1$ is a $K$-embedding such that $f(A_0) \subset A_1$, then $f\rest A_0$ is a $K_{dense}$-embedding from $(M_0)_{A_0}$ to $(M_1)_{A_1}$.
	\item Given $\mathcal{A}, \mathcal{B} \in K_{dense}$ and a $K_{dense}$-embedding $f: \mathcal{A} \to \mathcal{B}$, this lifts canonically to a $K$-embedding $\overline{f}:\overline{\mathcal{A}} \to \overline{\mathcal{B}}$.
	\end{itemize}
\end{enumerate}
\end{theorem}

{\bf Proof:} The proof proceeds similar to the first-order version, Theorem \ref{goal}.  In particular, many of the definitions of continuous structures from discrete approximations (such as getting the metric and relations from their approximations) did not use compactness.  Rathery, they only used uniform continuity to ensure that a completion existed, which will be guaranteed by the definition of $K_{dense}$ in this case.

Given continuous $\tau = \seq{F_i, R_j}_{i < n_f, j < n_r}$, define
$$\tau^+ := \seq{F_i^+, R_{R_j(\bz) \geq r}, R_{R_j(\bz) \leq r}}_{i < n_f, j< n_r, r \in \qp}$$
Given an MAEC $K$, we define the AEC $K_{dense}$ as follows:
\begin{itemize}
	\item $\tau(K_{dense}) = \tau^+$;
	\item Given an $\tau^+$ structure $\mathcal{A}$, we use the following procedure to determine membership in $K_{dense}$: define $D$ on $A \times A$ by
	\begin{eqnarray*}
	D(a, b) &:=& \sup\{ r \in \qp : \mathcal{A} \vDash R_{d(x, y) \geq r}[a, b]\}\\
	&=& \inf\{ r \in \qp: \mathcal{A} \vDash R_{d(x, y) \leq r}[a, b]\}
	\end{eqnarray*}
	The proof that this is a well-defined and is a metric proceeds exactly as in the previous case.  We can similarly define the relations $R_j$ on $A$ and complete the universe $(A, D)$ to $\overline{A}$.  We call the structure $\mathcal{A}$ \emph{completable} iff
	\begin{enumerate}
		\item for every $\ba \in \overline{A}$ and every Cauchy sequence $\seq{\ba^n \in A : n < \omega}$ converging to to $\ba$, $\lim_{n \to \infty}R_j(\ba^n)$ exists and its value is independent of the choice of the sequence; and similarly
		\item for every $\ba \in \overline{A}$ and every Cauchy sequence $\seq{\ba^n \in A : n < \omega}$ converging to to $\ba$, $\lim_{n \to \infty}F^+_i(\ba^n)$ exists and its value is independent of the choice of the sequence.
	\end{enumerate}
	If $\mathcal{A}$ is completable, then we define $\overline{\mathcal{A}}$ to be the $\tau^+$-structure where $F_i$ and $R_j$ are defined on $\overline{A}$ according to the independent value given above.  Finally, we say that $\mathcal{A} \in K_{dense}$ iff
	\begin{enumerate}
		\item $\mathcal{A}$ is completable; and
		\item $\overline{\mathcal{A}} \in K$.
	\end{enumerate}
	
	\item Given $\mathcal{A}, \mathcal{B} \in K_{dense}$, we say that $\mathcal{A} \prec_{dense} \mathcal{B}$ iff
	\begin{enumerate}
		\item $\mathcal{A} \subset_{\tau^+} \mathcal{B}$; and
		\item $\overline{\mathcal{A}} \prec_K \overline{\mathcal{B}}$.
	\end{enumerate}
\end{itemize}

Now that we have the definition of $K_{dense}$, we must show that it is in fact an AEC.  The verification of the axioms are routine.  We give the arguments for coherence and the chain axioms as templates.

For coherence, suppose that $\mathcal{A}, \mathcal{B}, \mathcal{C} \in K_{dense}$ such that $\mathcal{A} \prec_{dense} \mathcal{C}$; $\mathcal{B} \prec_{dense} \mathcal{C}$; and $\mathcal{A} \subset_{\tau^+} \mathcal{B}$.  Then, taking completions, we get that
$$\overline{\mathcal{A}} \prec_K \overline{\mathcal{C}}; \overline{\mathcal{B}} \prec_K \overline{\mathcal{C}}; \text{ and } \overline{\mathcal{A}} \subset_\tau \overline{\mathcal{B}}$$
By coherence in $K$, we then have that $\overline{\mathcal{A}} \prec_K \overline{\mathcal{B}}$.  Then, by definition, $\mathcal{A} \prec_{dense} \mathcal{B}$, as desired.

For the chain axioms, suppose that $\seq{\mathcal{A}_i \in K_{dense} : i < \alpha}$ is a continuous, $\prec_{dense}$-increasing chain such that, for all $i < \alpha$, $\mathcal{A}_i \prec_{dense} \mathcal{B} \in K_{dense}$.  Again, taking completions, we get that $\seq{\overline{\mathcal{A}}_i \in K : i < \alpha}$ is a continuous, $\prec_K$-increasing chain such that, for all $i < \alpha$, $\overline{\mathcal{A}}_i \prec_K \overline{\mathcal{B}} \in K$.  Then, by the union axioms for $K$, we have that $\overline{\cup_{i < \alpha} \overline{\mathcal{A}}_i} \in K$ and $\overline{\cup_{i < \alpha} \overline{\mathcal{A}}_i} \prec_K \overline{\mathcal{B}}$.  Note that the existence of $\overline{\cup_{i < \alpha} \overline{\mathcal{A}}_i}$ shows that $\cup_{i < \alpha} \mathcal{A}_i$ is completable and that an easy computation shows that $\overline{\cup_{i < \alpha} \overline{\mathcal{A}}_i} = \overline{\cup_{i < \alpha} \mathcal{A}_i}$.  Thus, $\cup_{i < \alpha} \mathcal{A}_i \in K_{dense}$ and, for all $j < \alpha$, 
$$\mathcal{A}_j \prec_{dense} \cup_{i < \alpha} \mathcal{A}_i \prec_{dense} \mathcal{B}$$

Once we have defined the maps and shown that $K_{dense}$ is an AEC, the rest of the proof proceeds exactly as in the continuous first-order case.  These proofs are sometimes simpler since $\tau^+$ only has relations for each relation of $\tau$, rather than each formula of $\tau$.\hfill\dag \\

We now turn to an application.  Both Hirvonen and Zambrano have proven versions of Shelah's Presentation Theorem for MAECs in their theses.  The more general is Zambrano's \cite[Theorem 1.2.7]{pedrothesis}:
\begin{theorem}
Let $K$ be a MAEC.  There is $\tau_1 \supset \tau$ and a continuous\footnote{But not \emph{uniformly} continuous.} $\tau_1$-theory $T_1$ and a set of $T_1$-types $\Gamma$ such that $K = PC(T_1, \Gamma, \tau)$.
\end{theorem}

An immediate corollary to our presentation theorem is a discrete presentation theorem.
\begin{cor}[Discrete Presentation Theorem for Metric AECs]
Let $K$ be a MAEC.  Then there is a (discrete) language $\tau_1$ of size $LS(K)$, an $\tau_1$-theory $T_1$, and a set of $T_1$-types $\Gamma$ such that $K = \{ \overline{M_1 \rest \tau(K)} : M_1 \vDash T_1 \te{ and omits }\Gamma\}$, where the completion is taken with respect to a canonically definable metric.
\end{cor}

{\bf Proof:} Apply Theorem \ref{Kdense} to represent $K$ as a discrete AEC $K_{dense}$ and then apply Shelah's Presentation Theorem from Shelah \cite{sh88}.\hfill\dag\\

Additionally, Zambrano asks \cite[Question 1.2.9]{pedrothesis} if there exists is a Hanf number for model existence in MAECs.  Using our presentation theorem, we can answer this questions in the affirmative.  Furthermore, the Hanf number for MAECs is the same as for AECs.

\begin{theorem} \label{maechanf-thm}
 If $K$ is a MAEC with $LS(K) = \kappa$  and has models  of size or density character cofinal in $\beth_{(2^\kappa)^+}$, then $K$ has models with density character arbitrarily large.
\end{theorem}

{\bf Proof:} For every, $\lambda < \beth_{(2^\kappa)^+}$, let $M_\lambda \in K$ have size $\geq \lambda$.  $M_\lambda$ is nicely dense in itself, so $(M_\lambda)_{|M_\lambda|} \in K_{dense}$ has size $\geq \lambda$.  By the definition of Hanf number for discrete AECs, this means that $K_{dense}$ has arbitrarily large models.  Taking completions, this means $K$ has arbitrarily large models.  The proof for density character is the same.\hfill \dag\\

Given this representation, we can determine basic structural properties of $K$ by looking at $K_{dense}$ and vice versa.  The above theorem already shows how to transfer arbitrary large models.  Other properties transfer similarly.

\begin{prop}
Suppose $K$ is an MAEC.  For $P$ being ``amalgamation,'' ``joint embedding,'' or ``no maximal models,'' $K$ has $P$ iff $K_{dense}$ has $P$.
\end{prop}

The statement of Theorem \ref{Kdense}.(4) is designed to make the proof of this proposition easy to see.

We now look at the notion of type in $K_{dense}$ that corresponds to Galois types in $K$.  Similar to the examination in Section \ref{type-sec}, we use the notion of sequence types.

\begin{defin}\
\begin{itemize}
	\item Given $\mathcal{A} \in K_{dense}$, $\seq{r_n : n < \omega}$ is a \emph{sequence Galois type over $\mathcal{A}$} iff
	\begin{enumerate}
		\item $r_0 \in gS^\ell(\mathcal{A})$; and
		\item $r_{n+1} \in gS^{\ell 2}(\mathcal{A})$ such that 
		\begin{enumerate}
			\item if $\ba \bb \vDash r_{n+1}$, then $d(\ba, \bb) \leq \frac{1}{2^n}$; and
			\item $r_{n+1}^{\{\ell, \dots, \ell 2 - 1\}} = r_n^\ell$, i. e. the first $\ell$ coordinates of $r_n$ are the same as the final $\ell$ coordinates of $r_{n+1}$.
		\end{enumerate}
	\end{enumerate}
	
	\item Given a sequence Galois type $\seq{r_n : n < \omega}$ over $\mathcal{A}$ and $\mathcal{A} \prec_{dense} \mathcal{B}$, we say that \emph{$\seq{\ba_n \in \mathcal{B} : n < \omega}$ realizes $\seq{r_n:n< \omega}$} iff
	\begin{enumerate}
		\item $\ba_0 \vDash r_0$; and
		\item $\ba_{n+1} \ba_n \vDash r_{n+1}$
	\end{enumerate}
	
	\item Given a sequence Galois type $\seq{r_n : n < \omega}$ over $\mathcal{A}$ and $\mathcal{A} \prec_{dense} \mathcal{B}$, we say that \emph{$\seq{\ba_n \in \mathcal{B} : n < \omega}$ weakly realizes $\seq{r_n:n< \omega}$} iff there is some $\mathcal{C}$ and $\seq{\bb_n \in \mathcal{C}:n<\omega}$ such that 
	\begin{enumerate}
		\item $\mathcal{B} \prec_{dense} \mathcal{C}$;
		\item $\seq{\bb_n : n < \omega}$ realizes $\seq{r_n : n < \omega}$; and
		\item $\lim_{n<\omega} \ba_n = \lim_{n<\omega} \bb_n$.
		
	\end{enumerate}
	
\end{itemize}
\end{defin}

Note that realizing a sequence Galois type is different than realizing each individual Galois type in the sequence separately.  However, we can always realize sequence Galois types given amalgamation.

\begin{lemma}\label{realize-seq-type}
Suppose that $K$ has amalgamation.  Given any sequence Galois type $\seq{r_n:n<\omega}$ over $\mathcal{A} \in K_{dense}$, there is $\mathcal{A} \prec_{dense} \mathcal{B} \in K_{dense}$ that contains a realization of $\seq{r_n:n< \omega}$.
\end{lemma}

{\bf Proof:}  By the definition of Galois type, we can write each $r_n$ as $gtp_{K_{dense}}(\ba^0_0/\mathcal{A}; \mathcal{B}_0)$ and $gtp_{K_{dense}}(\ba^{n+1}_{n+1} \ba^{n+1}_n/\mathcal{A}; \mathcal{B}_{n+1})$.  Using amalgamation in $K_{dense}$, which follows from amalgamation in $K$, we can construct increasing $\seq{\mathcal{C}_n:n<\omega}$ and increasing $f_n: \mathcal{B}_n \to_{\mathcal{A}} \mathcal{C}_n$ such that $f_0$ is the identity and $f_n(a^n_n) = f_{n+1}(a^{n+1}_n)$.  This second part is due to the second clause in the definition of sequence types.  Setting $\mathcal{C} = \cup_{n < \omega} \mathcal{C}_n$, we have that $\seq{f_n(a^n_n) : n < \omega}$ realizes $\seq{r_n : n < \omega}$. \hfill \dag\\

Unfortunately, we don't have the same tight connection between Galois types in $K$ and sequence Galois types in $K_{dense}$ as exists in Theorem \ref{limittypes}.  This is due to the fact that the Galois version of sequence types specifies a particular distance between consecutive members of the Cauchy sequence, rather than just specifying a bound on the distance.  It is possible that perturbations (as in Hirvonen and Hyttinnen \cite{maecpert}) might be used to restore this connection.  Instead, we have introduced the notion of weakly realizing a sequence Galois type because this is enough to prove a variant of Theorem \ref{Tdense-sat} in this context.

Following the L\"{o}wenheim-Skolem axiom, Galois saturation in the context of MAECs considers the density character of the domain of the type (see \cite[Definition 3.5]{maecs}).  This is important because we cannot consider Galois types over arbitrary sets as we could in the first-order case.

\begin{theorem}\label{Kdense-sat}
Let $M \in K$ and $\mathcal{A} \in K_{dense}$ such that $\overline{\mathcal{A}} = M$.  $M$ is $\kappa$-Galois saturated iff $\mathcal{A}$ is $\kappa$-weakly saturated for sequence Galois types.
\end{theorem}

{\bf Proof:} First, suppose that $M$ is $\kappa$ saturated.  Let $\mathcal{A}_0 \prec_{dense} \mathcal{A}$ be of size $< \kappa$ and let $\seq{r_n : n < \omega}$ be a sequence Galois type over $\mathcal{A}_0$.  By Lemma \ref{realize-seq-type}, there is some $\mathcal{B} {}_{dense} \succ \mathcal{A}_0$ and $\seq{\ba_n : n < \omega}$ that realizes $\seq{r_n : n < \omega}$.  After completing the members of $K_{dense}$, we have $\overline{\mathcal{A}_0} \prec M$ and $\overline{\mathcal{A}_0} \prec \overline{\mathcal{B}} \in K$.  Since $\seq{\ba_n \in \mathcal{B} : n < \omega}$ is a Cauchy sequence, there is some $\ba \in \overline{\mathcal{B}}$ such that $\lim_{n \to \infty} \ba_n= \ba$.  Since $M$ is $\kappa$-Galois saturated and $dc(\overline{\mathcal{A}_0}) < \kappa$, there is $\bb \in M$ that realizes $gtp_K(\ba/\overline{\mathcal{A}_0}; \overline{\mathcal{B}})$.  Since $\mathcal{A}$ is dense in $M$, there is some Cauchy sequence $\seq{\bb_n \in \mathcal{A} : n < \omega}$ that converges to $\bb$.  Then $\seq{\bb_n:n<\omega}$ weakly realizes $\seq{r_n : n < \omega}$.

Second, suppose that $\mathcal{A}$ is $\kappa$-weakly saturated for sequence Galois types and let $M_0 \prec M$ of density character $< \kappa$ and $r \in gS(M_0)$.  Let $A_0 \subset A$ be nicely dense in $M_0$.  Then $\mathcal{A}_0 := (M_0)_{A_0} \prec_{dense} \mathcal{A}$.  In $K$, we can write $r$ as $gtp_K(\ba/M_0; N)$.  Then, in $K_{dense}$, $\overline{r} := \seq{gtp_{K_{dense}}(\ba/ \mathcal{A}_0; N_{|N|}) : n < \omega}$ is a sequence Galois type over $\mathcal{A}_0$.  Since $\|\mathcal{A}_0\| < \kappa$, by $\mathcal{A}$'s weak saturation, there is some $\seq{\bb_n: n < \omega}$ that weakly realizes $\overline{r}$.  This means that we can find an extension $\mathcal{B} {}_{dense}\succ \mathcal{A}$ and $f:N_{|N|} \to_{\mathcal{A}_0} \mathcal{B}$ such that 
$$\lim_{n \to \infty} f(\ba) = \lim_{n \to \infty} \bb_n$$
Since $\bb_n \in \mathcal{A}$, this means that $f(\ba) \in M$.  Since $N_{|N|}$ is complete, the $K_{dense}$-embedding $f$ is in fact a $K$-embedding from $N$ into $\overline{\mathcal{B}}$ and fixes $\overline{\mathcal{A}_0} = M_0$.  Thus, we have that $f(\ba) \in M$ realizes $gtp_K(\ba/M_0; N)$. \hfill \dag\\

A crucial property in the study of MAECs is whether the natural notion of distance between Galois types defines a metric or not.  This property is called the Pertubation Property by Hirvonen and Hyttinen \cite{maecs} and the Continuity Type Property by Zambrano \cite{zamdtame}.  For ease, we assume that $K$ (equivalently, $K_{dense}$) has a monster model $\sea$.
\begin{defin}\
\begin{itemize}
	\item Given $M \in K$ and $p, q \in S(M)$, we define
		$$\bd(p, q) = \inf \{ d(a, b) : a\vDash p \text{ and }b \vDash q \}$$
	\item $K$ has the Pertubation Property (PP) iff, given any Cauchy sequence $\seq{b_n \in \sea : n < \omega}$ and $M \in K$, if, for all $n < m < \omega$, 
	$$gtp(b_n/M) = gtp(b_m/M)$$
	then, $gtp(\lim_{n \to \infty} b_n / M) = gtp(b_0/M)$.
\end{itemize}
\end{defin}

Although these properties might not initially seem related, a little work shows that $\bd$ is always a pseudometric and that it is a metric iff PP holds \cite{maecs}.  Then, similar to tameness from AECs, Zambrano \cite{zamdtame}.2.9 defines a notion of tameness in MAECs that satisfy PP.

\begin{defin}
$K$ is \emph{$\mu$-d-tame} iff for every $\epsilon > 0$, there is a $\delta > 0$ such that for every $M \in K$ of density character $\geq \mu$ and $p, q \in gS(M)$, if $\bd(p, q) \geq \epsilon$, then there is some $N \prec M$ of density character $\mu$ such that $\bd(p\rest N, q \rest N) \geq \delta_\epsilon$.
\end{defin}

Again these properties transfer to sequence types in $K_{dense}$.  We can define a pseudometric on sequence Galois types in $K_{dense}$ by
\begin{eqnarray*}
\bd_{dense}(\seq{r_n:n<\omega}, \seq{s_n:n<\omega}) &:=& \inf \{ \lim_{n\to \infty}d(a_n, b_n) : \seq{a_n : n < \omega} \text{ realizes } \\ 
& &\seq{r_n : n < \omega}, \seq{b_n:n<\omega} \text{ realizes }\seq{s_n : n<\omega}\}
\end{eqnarray*}

Then $\bd_{dense}$ is a metric iff $\bd$ is, and $\mu$-d-tameness transfers from $K$ to $K_{dense}$ in the obvious way.  This functor has recently been used in Boney and Zambrano \cite{metric-tamelc} to transfer various large cardinal results from Boney \cite{tamelc} to the MAEC setting.  Additionally, it can be used to transfer stability theoretic results both via characterization (as in Section \ref{type-sec}) and via independence relations.

\end{document}